\documentclass{elsart3-1}

\usepackage[applemac]{inputenc}
\usepackage{amssymb}
\usepackage[english,francais]{babel}

\newtheorem{theorem}{Theorem}[section]

\newtheorem{e-proposition}[theorem]{Proposition}

\newtheorem{e-definition}[theorem]{Definition\rm}

\newtheorem{theoreme}{Th\'eor\`eme}[section]

\renewcommand{\theequation}{\arabic{equation}}
\def\og{\leavevmode\raise.3ex\hbox{$\scriptscriptstyle\langle\!\langle$~}}
\def\fg{\leavevmode\raise.3ex\hbox{~$\!\scriptscriptstyle\,\rangle\!\rangle$}}

\def\R{\mathbb{R}}
\def\e{\varepsilon}
\def\Fe{F_\varepsilon}
\def\xe{x_\varepsilon}
\def\te{t_\varepsilon}
\def\we{w_\varepsilon}
\def\ue{u^\varepsilon}
\def\me{m_\varepsilon}
\def\tme{\tilde m_\varepsilon}
\def\ued{u^\varepsilon_\delta}
\def\ou{\overline u}
\def\oud{\overline u_\delta}
\def\uu{\underline u}
\def\xb{\bar x}
\def\tb{\bar t}
\def\SYM{\mathcal{S}(n)}
\def\YtoXandEPStoZERO{{\displaystyle{\mathop{\scriptstyle{(y,s)\to
(x,t)}}_{\varepsilon\to 0}}}} 
\def\limssup{\mathop{\rm lim\,sup\!^*\,}} \def\limiinf{\mathop{\rm 
lim\,inf_*\,}}

\begin{document}
\centerline{}
\begin{frontmatter}


\selectlanguage{english}
\title{A New Stability Result for Viscosity Solutions of Nonlinear Parabolic Equations with Weak Convergence in Time}

\selectlanguage{english}
\author[authorlabel1]{Guy Barles}
\ead{barles@lmpt.univ-tours.fr}
\address[authorlabel1]{Laboratoire de Math\'ematiques et Physique Th\'eorique (UMR CNRS 6083), Fédération Denis Poisson, Universit\'e de Tours,
Parc de Grandmont, 37200 Tours, France}


\medskip
\begin{center}
{\small Received *****; accepted after revision +++++\\
Presented by Pierre-Louis Lions}
\end{center}

\begin{abstract}
\selectlanguage{english}
We present a new stability result for viscosity solutions of fully nonlinear parabolic equations which allows to pass to the limit when one has only weak convergence in time of the nonlinearities.
{\it To cite this article: G. Barles, C. R. Acad. Sci. Paris, Ser. I 340 (2005).}

\vskip 0.5\baselineskip

\selectlanguage{francais}
\noindent{\bf R\'esum\'e} \vskip 0.5\baselineskip \noindent
{\bf Un nouveau résultat de stabilité pour les solutions de viscosité d'équations paraboliques non-linéaires avec convergence faible en temps}
Nous obtenons un nouveau résultat de stabilité pour les solutions de viscosité d'équations fortement non linéaires paraboliques dans le cas où l'on n'a qu'une convergence faible en temps pour les non-linéarités.{\it Pour citer cet article~: G. Barles, C. R. Acad. Sci.
Paris, Ser. I 340 (2005).}

\end{abstract}
\end{frontmatter}

\selectlanguage{francais}
\section*{Version fran\c{c}aise abr\'eg\'ee}
Le but de cette Note est de décrire un nouveau résultat de stabilité pour les solutions de viscosité d'équations fortement non-linéaires paraboliques avec des dépendances mesurables en temps ($L^1$ ou $L^\infty$); l'originalité de ce résultat est de prendre en compte des convergences faibles en temps pour les non-linéarités au lieu de l'hypothèse classique de convergence forte dans $L^1$.

Pour \^etre plus précis, nous considérons une suite $(\ue)_\e$ de solutions de~:
\begin{equation}\label{Fee}
 \ue_t + \Fe (x,t, \ue, D\ue, D^2 \ue) = 0 \quad \hbox{dans  }\Omega \times (0,T),
\end{equation}
où $\Omega$ est un ouvert de $\R^n$ et $T>0$. Les solutions $\ue$ sont scalaires, 
$ \ue_t $, $D\ue$, $D^2\ue$ désignent respectivement la dérivée en temps, le gradient et la matrice hessienne de $\ue$.
Les non-linéarités $\Fe (x,t,r,p,X) $ sont des fonctions à valeurs réelles qui sont définies pour presque tout $t\in (0,T)$ et pour tout $(x,r,p,X)\in \Omega \times \R\times \R^n \times \SYM$ où $\SYM$ est l'espace des matrices symétriques $n\times n$. On suppose que ces fonctions sont continues en $(x,r,p,X)$ pour presque tout $t$, que $t\mapsto \Fe (x,t,r,p,X)\in L^{1}(0,T)$ pour tout $(x,r,p,X)$ et qu'elles sont elliptiques dégénérées au sens suivant~:

\smallskip
$ \Fe (x,t,u,p,X) \geq \Fe (x,t,u,p,Y)$ si $X\leq Y$, pour presque tout $t \in (0,T)$ et pour tous $x\in \Omega$, $u\in \R$, $p\in \R^n$ et pour toutes matrices $X,Y$.

\smallskip
Avant de présenter le résultat principal, nous rappelons que la notion de solution de viscosité pour des équations avec une dépendance $L^1$ en temps a été introduite par H. Ishii \cite{I} pour les équations de Hamilton-Jacobi du premier ordre. Dans ce cadre, la définition n'est pas une généralisation triviale du cas classique où l'équation est continue en temps; malgré cette difficulté, H. Ishii obtenait, dans \cite{I}, des résultats d'existence, d'unicité et de stabilité tout à fait analogues à ceux du cas classique et sous des hypothèses naturelles. Nous renvoyons aussi à P.L. Lions et B. Perthame \cite{LP} pour une approche différente mais équivalente à celle d'Ishii. Pour les équations du second ordre, les premiers travaux sont dus à D.~Nunziante \cite{DN1,DN2}; le lecteur pourra aussi consulter les travaux plus récents de M. Bourgoing\cite{MB1,MB2}. 

Pour formuler les hypothèses principales, nous introduisons la propriété suivante~: nous dirons que la fonction $H : \Omega \times (0,T) \times \R \times \R^n \times \SYM \to \R$ satisfait la propriété M(K) pour un compact $K$ de $\Omega \times \R\times \R^n \times \SYM$ s'il existe un module de continuité $m=m(K) : (0,T) \times \R^{+}\rightarrow \R^{+}$ tel que $t\mapsto m (t,r)\in L^1(0,T)$ pour tout $ r\geq 0$, $m (t,r)$ est croissant par rapport à $r$, $m (.,r) \to 0$ in $L^1(0,T)$ quand $r \to 0$ et~:

\smallskip
\centerline{$\displaystyle
|H (x_1,t,r_1,p_1,X_1)- H (x_2,t,r_2,p_2,X_2)| \leq m (t, |x_1-x_2| + |u_1-u_2| + |p_1-p_2| +|X_1-X_2| ),$}

\smallskip
\noindent pour presque tout $t \in (0,T)$ et pour tous $(x_1,r_1,p_1,X_1), (x_2,r_2,p_2,X_2) \in K$. De plus, nous dirons que $H$ satisfait la propriété (M) s'il satisfait M(K) pour tout compact $K$ de $\Omega \times \R\times \R^n \times \SYM$.

Les hypothèses sur les $(\Fe)_\varepsilon$ sont les suivantes~:\\
{\bf (F1)} Pour tout $\varepsilon >0$, $\Fe$ satisfait la propriété (M) avec des modules de continuité $\me=\me(K)$ qui satisfont $|| \me (.,r) ||_{L^1(0,T)} \to 0$ quand $r \to 0$ uniformément en $\varepsilon$, pour tout compact $K$. \\
{\bf (F2)} Il existe une fonction $F$ satisfaisant la propriété (M) telle que, pour tout $(x,r,p,X)$~:\\
\centerline{$ \int_0^t \Fe (x,s,r,p,X)ds \to \int_0^t  F (x,s,r,p,X)ds \quad\hbox{localement uniformément dans $(0,T)$}.$}

\smallskip
Le résultat principal est le suivant~:
\begin{theoreme}\label{mrf} On suppose que $(\Fe)_\varepsilon$ est une suite de fonctions elliptiques dégénérées qui satisfont {\bf (F1)}-{\bf (F2)}. Si $(\ue)_\varepsilon$ est une suite de sous-solutions (resp. sursolutions) uniformément localement bornées de (\ref{Fee}), alors $\ou:=\limssup \!\ue$ défini par (\ref{slr}) (resp. $\uu:=\limiinf \!\ue$) est une sous-solution (resp. sursolution)~de~:
$$ w_t + F (x,t, w, Dw, D^2 w) = 0 \quad \hbox{in  }\Omega \times (0,T).$$
\end{theoreme}
 
Ce résultat s'étend aux équations singulières de type géométriques, intervenant typiquement dans l'approche par ligne de niveaux pour définir des mouvements d'interfaces avec des vitesses dépendant de la courbure: cette généralisation ne nécessite que des modifications naturelles des hypothèses {\bf (F1)}-{\bf (F2)}. Il se généralise aussi, sans difficulté, à tous les types classiques de conditions aux limites à condition que $u_t$ n'intervienne pas dans la condition aux limites.

\selectlanguage{english}

\section{Introduction~: the main result}
\label{}

The aim of this Note is to provide a new stability result for viscosity solutions of fully nonlinear parabolic equations with measurable dependences in time ($L^1$ or $L^\infty$); the originality of this new result is to handle the case of weak convergence in time of the equations, instead of the classical strong $L^1$ convergence.

In order to be more specific, we consider a sequence of fully nonlinear parabolic equations
\begin{equation}\label{Fe}
 \ue_t + \Fe (x,t, \ue, D\ue, D^2 \ue) = 0 \quad \hbox{in  }\Omega \times (0,T),
\end{equation}
where $\Omega$ is an open subset of $\R^n$, $T>0$ and the solutions $\ue$ are scalar; here and below 
$ \ue_t $ denotes the derivative of $\ue$ with respect to
$t$ and $D\ue$, $D^2\ue$ denote respectively the gradient and the Hessian
matrix of $\ue$ with respect to the space variable $x$. 
The nonlinearities $\Fe (x,t,r,p,X) $ are real-valued functions, defined for almost every $t\in (0,T)$ and for every $(x,r,p,X)\in \Omega \times \R\times \R^n \times \SYM$, where $\SYM$ denotes the space of $n\times n$ symmetric matrices. These functions are assumed to be continuous with respect to $(x,r,p,X)$ for almost every $t$,
$t\mapsto \Fe (x,t,r,p,X)\in L^{1}(0,T)$ for every $(x,r,p,X)$ and they are degenerate elliptic, i.e

\smallskip
$ \Fe (x,t,u,p,X) \geq \Fe (x,t,u,p,Y)$ if $X\leq Y$, for almost every $t \in (0,T)$ and for any $x\in \Omega$, $u\in \R$, $p\in \R^n$ and any $n \times n$-symmetric matrices $X,Y$.

\smallskip
Before presenting our main result, we recall that the notion of viscosity solutions for such equations with a measurable dependence in time was first introduced by H. Ishii \cite{I} for first-order Hamilton-Jacobi Equations; in this framework, the definition is not a straightforward extension of the classical case (the case of equations with continuous dependence in time) but, in \cite{I}, H. Ishii obtained natural generalizations of the classical existence, uniqueness and stability results. A different but equivalent approach was then provided by P.L. Lions and B. Perthame \cite{LP}. The extension to the case of second-order equations was first done by D.~Nunziante \cite{DN1,DN2} (see also the more recent works of M. Bourgoing\cite{MB1,MB2}). Finally, for a wider presentation of the theory, we refer the reader to the User's guide of M.G. Crandall, H. Ishii and P.L Lions\cite{CIL} and the books of M.~Bardi, I.~Capuzzo-Dolcetta\cite{BC} and W. H.~Fleming and H. M.~Soner\cite{FS} .

In order to formulate the needed assumptions on the $\Fe$'s, we introduce the following property~: we say that a function $H : \Omega \times (0,T) \times \R \times \R^n \times \SYM \to \R$ satisfies the property M(K) on a compact subset $K$ of $\Omega \times \R\times \R^n \times \SYM$ if there exists a modulus $m=m(K) : (0,T) \times \R^{+}\rightarrow \R^{+}$ such that $t\mapsto m (t,r)\in L^1(0,T)$ for all $ r\geq 0$, $m (t,r)$ is non-decreasing in $r$, $m (.,r) \to 0$ in $L^1(0,T)$ as $r \to 0$,  and
$$
 |H (x_1,t,r_1,p_1,X_1)- H (x_2,t,r_2,p_2,X_2)| \leq m (t, |x_1-x_2| + |u_1-u_2| + |p_1-p_2| +|X_1-X_2| ),
$$
for almost every $t\in (0,T)$ and for any $(x_1,r_1,p_1,X_1), (x_2,r_2,p_2,X_2) \in K$. Moreover, we say that $H$ satisfies the property (M) if it satisfies M(K) for any compact subset $K$ of $\Omega \times \R\times \R^n \times \SYM$.

The assumptions on the sequence $(\Fe)_\varepsilon$ are the following\\
{\bf (F1)} For any $\varepsilon >0$, $\Fe$ satisfies the property (M) for some modulus $\me=\me(K)$ such that $|| \me (.,r) ||_{L^1(0,T)} \to 0$ as $r \to 0$ uniformly with respect to $\varepsilon$, for any compact subset $K$. \\
{\bf (F2)} There exists a function $F$ satisfying the property (M) such that, for any $(x,r,p,X)$, 
$$ \int_0^t \Fe (x,s,r,p,X)ds \to \int_0^t  F (x,s,r,p,X)ds \quad\hbox{locally uniformly in $(0,T)$.}$$

Finally we recall the following classical definition of half-relaxed limits~: if $(\ue)_\varepsilon$ is a sequence of locally uniformly bounded functions then we set, for $(x,t) \in \Omega \times (0,T)$
\begin{equation}\label{slr}
 \limssup \ue(x,t)=\limsup_\YtoXandEPStoZERO\,\ue (y,s) \quad\hbox{and}\quad
\limiinf \ue(x,t) =\liminf_\YtoXandEPStoZERO\,\ue (y,s)\; .
\end{equation}

\medskip
The result is the following.

\begin{theorem}\label{mr} Assume that $(\Fe)_\varepsilon$ is a sequence of degenerate elliptic functions satisfying {\bf (F1)}-{\bf (F2)}. If $(\ue)_\varepsilon$ is a sequence of locally uniformly bounded subsolutions (resp. supersolutions) of (\ref{Fe}), then $\ou:=\limssup \ue$ (resp. $\uu:=\limiinf \ue$) is a subsolution (resp. supersolution) of
\begin{equation}\label{F}
 w_t + F (x,t, w, Dw, D^2 w) = 0 \quad \hbox{in  }\Omega \times (0,T).
\end{equation}

\end{theorem}

\medskip
\noindent{\bf Proof :} we only prove the result for the subsolution case, the supersolution case being proved analogously.

We use Ishii's definition of viscosity solutions for equations with a $L^1$ dependence in time : again we refer to D. Nunziante\cite{DN1,DN2} and M. Bourgoing\cite{MB1,MB2} for a complete presentation of the theory.

Let $(x,t)$ be a strict local maximum point of $\ou - \phi + \int_0^t b(s)ds$ where $\phi$ is a smooth test-function and $b$ is a $L^1$-function  such that there exists a continuous function $G$ satisfying $ b(s) + G(y,s,v,q,M) \leq F (y,s,v,q,M)$ in a neighborhood of $(x,t, \ou (x,t), D\phi (x,t) , D^2\phi (x,t))$ a.e. in $s$ and for all $(y,v,q,M)$. We have to prove that
$$ \phi_t (x,t) + G(x,t, \ou (x,t), D\phi (x,t) , D^2\phi (x,t)) \leq 0.$$

To do so, we first pick some small $\delta >0$ and we consider $m, \me $ given by {\bf (F1)} for a large compact $K$. We introduce a new sequence $(\ued)_\varepsilon$ defined by
$$ \ued (x,t) := \ue (x,t) - \int_0^t \, [\me (s,\delta) + m(s,\delta)]ds\; .$$
We denote by $\oud = \limssup  \ued$. By the properties of $\me$ and $m$, we have $\uu - o_\delta (1) \leq \oud \leq \uu $ in $\R^n \times (0,T)$.

Using this last property, by classical results, since $(x,t)$ is a \emph{strict} local maximum point of $\ou - \phi + \int_0^t b(s)ds$, if $\delta$ is sufficiently small, there exists a local maximum point $(\xb,\tb)$ of $\oud - \phi + \int_0^t b(s)ds$ near $(x,t)$. We have dropped the dependence of $(\xb,\tb)$ in $\delta$ for the sake of simplicity of notations but clearly $(\xb,\tb) \to (x,t)$ when $\delta \to 0$. Moreover, subtracting if necessary a term like $(s-\tb)^2 + |y-\xb|^4$, we may assume that $(\xb,\tb)$ is a strict local maximum point as well.

Then we fix $\delta$ and we now consider the functions
$$ \chi_\varepsilon(y,s) := \ued (y,s)  - \phi (y,s) + \int_0^s b(\tau)d\tau +  \int_0^s \psi_\varepsilon (\tau )d\tau ,$$
where
$$ \psi_\varepsilon (\tau) := \Fe (\xb,\tau, \oud (\xb,\tb), D\phi (\xb,\tb) , D^2\phi (\xb,\tb)) -F (\xb,\tau, \oud (\xb,\tb), D\phi (\xb,\tb) , D^2\phi (\xb,\tb)) .$$
By assumption {\bf (F2)}, $ \int_0^s \psi_\varepsilon(s) ds$ converges uniformly to $0$ and, passing if necessary to a subsequence, there exist local maximum points $(\xe,\te)$ of $\chi_\varepsilon$ such that, as $\varepsilon \to 0$, $(\xe,\te) \to (\xb,\tb)$ and  $\ue (\xe,\te) \to \oud (\xb,\tb)$. 

Next, we remark that $\ued$ is a subsolution of
$$ (\ued)_t + \Fe (x,t, \ued, D\ued, D^2 \ued) \leq -\me (t,\delta) - m(t,\delta) \quad \hbox{in  }\Omega \times (0,T),$$
and in order to apply the definition, we have to look at 
$$Q(y,s,v,q,M):= b(s) + \psi_\varepsilon (s) + G(y,s,v,q,M).$$
By the property of $b$ and $G$, we have
\begin{eqnarray*}
Q(y,s,v,q,M) & \leq & \psi_\varepsilon (s)  + F(y,s,v,q,M) \\
& \leq & \Fe (y,s,v,q,M) + [\Fe (\xb,s, \oud (\xb,\tb), D\phi (\xb,\tb) , D^2\phi (\xb,\tb)) - \Fe (y,s,v,q,M)] + \\
& & \left[F(y,s,v,q,M) -F (\xb,s, \oud (\xb,\tb), D\phi (\xb,\tb) , D^2\phi (\xb,\tb)) \right] \; .
\end{eqnarray*}
But, by {\bf (F1)} and the fact that $$(\xe, \ued (\xe,\te), D\phi (\xe,\te) , D^2\phi (\xe,\te))\to (\xb, \oud (\xb,\tb), D\phi (\xb,\tb) , D^2\phi (\xb,\tb))$$ when $\varepsilon \to 0$, the right hand-side of this inequality is less than $\Fe (y,s, v, q , M) + \me (s,\delta) + m (s,\delta)$ in a neighborhood of $(\xe, \te, \ued (\xe,\te), D\phi (\xe,\te) , D^2\phi (\xe,\te))$ and we can apply the definition of viscosity subsolution for $\ued$ which yields
$$ \phi_t (\xe,\te) + G(\xe, \te, \ued (\xe,\te), D\phi (\xe,\te) , D^2\phi (\xe,\te))\leq 0\; .$$
In order to conclude, we first let $\varepsilon \to 0$ using the continuity of $G$ and the properties of $\xe, \te, \ued (\xe,\te)$. Finally we let $\delta \to 0$ using that, by the same arguments, $(\xb,\tb)\to (x,t)$ and  $\oud (\xb,\tb) \to \ou (x,t)$. And the proof is complete.

\section{On various possible extensions}

We briefly present in this section some easy extensions of the above main results in several different frameworks. 

We begin with the case of \emph{``singular equations''} arising typically in the so-called level-sets approach when considering the motion of hypersurfaces with curvature dependent velocities. On a technical point of view, this means to take in account locally bounded functions $\Fe(x,t,p,M), F(x,t,p,M)$ independent of $u$, which are continuous in $(x, p, M)$ for $p\neq 0$ for almost every $t$ and with a possible singularity for $p=0$. We refer to \cite{MB1,MB2} for results for such equations in the $L^1$-framework.

In order to extend Theorem~\ref{mr} to this case, we have to replace {\bf (F1)}-{\bf (F2)} by changing (M) in (M-s) i.e. $\Fe , F$ satisfy M(K) for any compact subset of $\Omega \times (0,T) \times \R^n- \{0\} \times \SYM$ and with modulus $\me$ satisfying the same properties as in {\bf (F1)}. We have also to add\\
{\bf (F3)} There exists a neighborhood $V$ of $(0,0)$ in $\R^n \times \SYM$ such that, for any $\varepsilon >0$ and for any compact subset $\tilde K$ of $\R^n$, there exists a modulus $\tme=\tme(K)$ satisfying the same properties as $\me$
 such that
$$|\Fe (x,t,p,M)| \leq \tme (t, |p| + |M|) \quad\hbox{a.e. for $t\in (0,T)$, for any $x \in \tilde K$ and $(p,M) \in V$.}$$
And $F$ satisfies an analogous property. Finally, for {\bf (F2)}, the convergence property has to hold for all $(x,p,M)$, $p\neq 0$.

The proof of this case follows along the lines of the above proof if $D \phi (x,t)\neq 0$ while, if $D \phi (x,t)=0$, one has to use the result of Ch.~Georgelin and the author~\cite{BG} saying that one may assume w.l.o.g that $D^2 \phi (x,t)=0$; this allows to use {\bf (F3)}.

A second extension is the case of various \emph{boundary conditions} (in the viscosity sense)~: the case of Dirichlet conditions is a straightforward adaptation while, for possibly non linear Neumann boundary conditions, one has to avoid boundary conditions depending on $u_t$.

\section{A discussion of the assumptions and of the result on a typical example}

We discuss here the following example motivated by the study of stochastic pdes
\begin{equation}\label{aspde}
\ue_t + F(x,t, \ue , D\ue, D^2 \ue) + \dot \we (t) H(D\ue) = 0 \quad \hbox{in  }\R^n \times (0,T)\; ,
\end{equation}
where $F, H$ are, at least, continuous functions, $F$ being degenerate elliptic and the $\we$ are $C^1$-functions.

Such approximate problems were introduced by P.L. Lions and P.E. Souganidis \cite{LS1,LS2,LS3} to prove the existence of solutions of stochastic pdes and to obtain various properties. Typically $\we$ may be a smooth, pathwise approximation of the brownian motion. Of course, the key question is to understand the behavior of $\ue$ under the (apparently very weak) assumption that $\we$ converges uniformly on $[0,T]$ to some function~$w$. 

The first remark is that, under this assumption, the convergence property in {\bf (F2)} holds whatever $\we$ and $w$ are. Here we restrict ourselves to the $L^1-L^\infty$ framework which means that we assume that $\dot w \in L^1 (0,T)$. Even with this additional assumption, one may emphasize three different levels of difficulty.

First, if $||\dot \we||_{L^1}$ is uniformly bounded, then {\bf (F1)}-{\bf (F2)} hold and Theorem~\ref{mr} readily applies.
 
On the contrary, if $||\dot \we||_{L^1}$ is not uniformly bounded, {\bf (F1)} does not hold. But if we assume instead that, for any smooth initial datas $\psi$, the equation $u_t + H(Du) = 0$ has a smooth solution $S_H (t)\psi$ for small (positive and negative) time, depending smoothly in $\psi$, then the above strategy of proof applies, changing the terms ``$\phi (y,s) +  \int_0^s \psi_\varepsilon (\tau )d\tau$'' into $S_H (\we (s)- w (s))\phi(\cdot,s)$. Formally, it is easy to see that this ``exact'' term allows, in particular, to drop the term ``$\int_0^t \, [\me (s,\delta) + m(s,\delta)]ds$'' which is there to control the errors (the complete rigourous proof follows along the lines of the formal proof if $w$ is $C^1$ and just requires straightforward approximation arguments if $\dot w$ is just in $L^1$). Therefore, in this second case, we can also pass to the limit. It is worth pointing out that introducing the semi-group $S_H$ in such a way is precisely the type of arguments which is used by P.L. Lions and P.E. Souganidis\cite{LS1,LS3} for stochastic pdes.

Finally, if we are not in one of the two above cases, both (relatively close) stategies fail and this is not surprising. Indeed, in the case when $F\equiv 0$ and $w\equiv 0$, the question reduces to: does the limit of the $\ue$'s solve $u_t=0$? And the answer given by P.L. Lions and P.E. Souganidis\cite{LS2,LS3} is non trivial and not so natural~: it is true \emph{if and only if} $H$ is the difference of two convex functions, a property that the above strategies of proof cannot see.


\end{document}